\documentclass[15]{article}
\usepackage{amsmath}
\usepackage{graphicx}
\usepackage{amssymb}
\usepackage{amsthm}
\usepackage{float}

\begin{document}
\title{Chen and Casorati curvature inequalities for the submanifolds of quaternionic Kaehler manifolds endowed with Ricci quarter-symmetric metric connection.}
\author{Umair Ali Wani, Mehraj Ahmad Lone\\
National Institute Of Techonolgy, Srinagar,Kashmir
\date{}\\
Email:umair\_26btech16@nitsri.ac.in}
\maketitle
\textbf{Abstract:} In this paper, authors have established Chen's inequalties for the submanifolds of quaternionic Kaehler manifolds characterized by Ricci quarter-symmetric metric connection. Other than these inequalities, we have also derived generalized normalized Casorati curvature inequalities for the same submanifolds.
\section*{Introduction:} 
To provide the answers to Chern's open problem concerning the existence of minimal immersions into a Euclidean space of arbitrary dimension, Chen established the famous invariant $'\delta'$ $[1]$ $[2]$. The equation involving this invariant $'\delta'$ for the Riemannian manifold, M is given as:\\
\begin{equation}
\delta_M= \tau(p)- inf (K)(p).
\end{equation}
Here, $'\tau (p)'$ is the scalar curvature of M and $'K(p)'$ is the sectional curvature.\\
The inequality for the same $'\delta'$ for any Riemannian manifold, M is given as [2]:
\begin{equation}
\delta_M \leq \frac{n^2 (n-1)}{2(n-1)} \|H\|^2 +\frac{1}{2}(n+1)(n-2)c,
\end{equation}
where $'c'$ is the constant sectional curvature and  $'H'$ is the mean curvature vector for this n-dimensional manifold.\\
 He also derived an inequality which establishes a relationship between the extrinsic (K and $\tau$) and intrinsic invariants (squared mean curvature, H) of the submanifold. The same inequality is given as: $[3]$
 \begin{equation}
H^2 \geq \frac{4}{n^2} \{Ric(X)-(n-1)c\},
\end{equation}
where $Ric(X)$ is Ricci curvature of M at X.
 
With the basic insight provided by Chen, these inequalities were further extended to different manifolds under various metric and non-metric connections. Chen and Dillen studied the inequalities for Lagrangian submanifolds in complex space forms $[4]$. For Sasakian manifolds, the inequalities were studied by Mihai $[5]$. Arslan et al. studied the inequalities in contact manifolds $[6]$. The study on Lorentzian manifolds was done by Gulbahar et al. $[7]$. Mihai and Ozgur $[8]$ studied Chen inequalities on Riemannian manifolds with semi-symmetric metric connection while Zhang et al. $[9]$ $[10]$ studied geometric inequalities for the submanifolds of Riemannian manifolds of quasi-constant curvature with semi-symmetric metric connection. Inequalities for the submanifolds of real space forms with Ricci quarter-symmetric metric connection were studied by Nergis and Halil $[11]$. The inequalities of submanifolds in quaternionic space forms were studied by Yoon $[12]$.\\
Another solution to the problem of establishing relationship between the main extrinsic and intrinsic invariants is given by inequalities involving Casorati Curvatures. Casorati Curvature was originally established for the surfaces in Euclidean 3-space $[13]$ as a normalized sum of squared principal curvatures. It was extended to the general case of a submanifold of a Riemannian manifold as the normalized square of the length of second fundamental form of the submanifold by Decu et al.$[14]$ $[15]$. The inequalities involving this curvature has been studied for various submanifolds. Vilcu  $[16]$ obtained an optimal inequality for Casorati Curvature in Lagrangian submanifolds in complex space forms. Brubaker and Suceava  $[17]$ obatained the geometrical interpretation of Cauchy-Schwarz inequality in terms of Casorati Curvature. Lee et al. $[18]$ $[19]$ established the inequalities for Casorati curvature of submanifolds in generalized space forms endowed with semi-symmetric metric connection and Kenmotso space forms. The study of Casorati Curvature on holomorphic statistical manifolds of constant holomorphic curvature was done by Decu et al. $[20]$. Similarly, the inequalities on  slant submanifolds in metallic Riemannian space forms were established by Chaudhary and Blaga $[21]$. \\
 The concept of Ricci Quarter-symmetric metric connection was established by Golab$[22]$ and further developed by Mishra and Pandey $[23]$ and Kamilya and De $[24]$.\\
In this paper, the author has derived the optimal inequalities for the submanifolds of quaternionic Kaehler manifolds endowed with Ricci quarter-symmetric metric connection.\\\\
\section*{Preliminairies:}
Let $M^{'}$ be the $4m$-dimensional Riemannian manifold. This $M^{'}$ is called the Quaternionic Kaehler manifold if there exists 3-dimensional vector space of tensors of type $(1,1)$ with the local basis of almost Hermitian structure $\psi_1$,$\psi_2$ and $\psi_3$ such that:
\begin{equation}
\psi^2 _i=-I
\end{equation} 
\begin{equation}
\psi_{i}\psi_{i+1}=\psi_{i+2}=-\psi_{i}\psi_{i+1} 
\end{equation}
Where I is the identity transformation of tangent space $T_p X$ at point p for the vector $X$ of $M^{'}$.\\
And
\begin{equation}
\nabla^{`}_X \psi_i=q_{i+2}(X)\psi_{i+1}-q_{i+1}(X)\psi_{i+2}
\end{equation}
with $q_{i}$ as the local one form.\\
If X is any unit vector on $M^{'}$, then X, $\psi_1 (X)$, $\psi_2 (X)$ and $\psi_3 (X)$ form the orthonormal frame on $M^`$. Let the plane spanned by these vectors be $Q(X)$. If X, Y are two unit vectors in $M^{'}$ and let the plane spanned by these vectors be denoted by $\pi(X,Y)$, then this plane is called the  Quaternionic plane in $Q(X)$. The sectional curvature of this plane is called as quaternionic sectional curvature of $\pi$. If this sectional curvature is constant say, $4c$, then the quaternionic Kaehler manifold is called as quaternionic space form. Quaternionic Kaehler manifold is called quaternionic space form if and only if its curvature tensor follows the equation as: $[25]$
 
\begin{multline}
R^{*}(X,Y;Z)=c\{g(Y,Z)X-g(X,Z)Y + \sum_{i=1} ^{3} (g(\psi_i Y,Z)\psi_i X- g(\psi_i X,Z)\psi_i Y\\ -2g(\psi_i X,Y)\psi_i Z)\},
\end{multline}
where X,Y and Z are tangent vectors and g is the Riemannian metric of $M^{'}$.\\\\
For this manifold $M^{'}$, if $\nabla ^{"}$ is its linear connection then this connection is known as Ricci quarter symmetric metric connection if the torsion tensor $T^{'}$ satisfies the following equation:

\begin{equation}
T^{'} = \eta(Y)LX -\eta(X)LY,
\end{equation}
 where $\eta$ is a $1$-form and $L$ is a $(1,1)$ Ricci tensor defined by:
 \begin{equation}
 g(LX,Y)=S(X,Y),
 \end{equation}
 with S as Ricci tensor of $M^{'}$.\\
 With $\nabla ^{*}$ as Levi-Civita connection with respect to metric $g$, the Ricci quarter symmetric metric connection is given by $[22]$:
 \begin{equation}
 \nabla^{"}_X Y = \nabla^{*}_X Y + \eta(Y)LX -S(X,Y)P,
 \end{equation}
 where X and Y are vectors in $M^{'}$ with the unit  vector field P defined as:
 \begin{equation}
 g(P,X)=\eta(X).
 \end{equation}
 Let $N^{'}$ be an n-dimensional submanifold immersed in Riemannian manifold $M^{'}$, then there exist Ricci quarter-symmetric metric connection, $\nabla $ and Levi-Civita metric connection, $\nabla^{'}$ which are induced on this submanifold $N^{'}$. Let $R^{"}$ and $R^{*}$ be curvature tensors associated with $\nabla^{"}$ and $\nabla^{*}$ respectively Similarly, let $R$ and $R^{'}$ be curvature tensors associated with $\nabla$ and $\nabla^{'}$ respectively. $R^{*}$ has been already explained in equation $(7)$.\\\\
 For this manifold $M^{'}$, the curvature tensor with respect to Ricci quarter-symmetric metric $\nabla ^{"}$ connection is given as:
 \begin{multline}
  R^{"}(X,Y;Z) =  R^{*}(X,Y;Z)- M(Y,Z)LX + M(X,Z)LY -S(Y,Z)QX+S(X,Z)QY \\+ \eta(z)[(\nabla^{*} _X L)Y-(\nabla^{*} _Y L)X]-[(\nabla^{*}_X S)(Y,Z)-(\nabla^{*}_Y S)(X,Z)]P,
\end{multline}
where X,Y and Z are tangent vectors in $M^{'}$ and M is a tensor field of type $(0,2)$ defined by;
\begin{multline}
  M(X,Y)=g(QX,Y)=(\nabla^{*}_X \eta)Y- \eta(Y) \eta(LX)+ \frac{1}{2} \eta(P)S(X,Y).
\end{multline}
In the above equation, Q is a tensor field of type $(2,1)$ defined by,
  \begin{equation}
 QX= \nabla^{*}_X P- (\eta)(LX)P + \frac{1}{2} \eta(P)LX.
  \end{equation} 
  If $N^{'}$ is assumed to be Einstein manifold, then the Ricci curvature S for the linear connection $\nabla ^{*}$ is given as: $[11]$
  \begin{equation}
  S(X,Y)=\frac{\tau^{'}}{n}g(X,Y),
  \end{equation}
 with $\tau^{'}$ as the scalar curvature.\\
 Let $\pi$ be any two dimensional plane of $T_p N^{'}$ generated by $\{e_i, e_j\}$ at point $p$, then the sectional curvature of the section is given by:
\begin{equation}
K_{ij} =\frac{g(R(e_i,e_j;e_j,e_i)}{g(e_i,e_i)g(e_j,e_j)-g^{2}(e_i,e_j)}, 
\end{equation}
where $e_i$, $e_j$ are orthonormal vectors in $T_p N^{'}$.\\
 If $(e_1,...,e_n)$ is the orthonormal basis of $T_p N^{'}$, then at any point p, scalar curvature $\tau$ is given as:
 \begin{equation}
 \tau(p) = \sum_{1\leq i <j \leq n} K_{ij}.
\end{equation} 
The normalized scalar curvature $\rho$ for $N^{'}$ is given as:
\begin{equation}
\rho =\frac{2\tau}{n(n-1)}.
\end{equation}  
 Using the equation $(15)$ and the fact that connection is metric connection in equation $(10)$, we get:
 \begin{multline}
 R^{"}(X,Y;Z) =  R^{*}(X,Y;Z)-\frac{\tau^{'}}{n}\{ M(Y,Z)X - M(X,Z)Y\\ +g(Y,Z)QX-g(X,Z)QY \}.
 \end{multline}
 Contracting the above equation with respect to X, we get:
 \begin{equation}
 S^{"}(Y,Z) =  \frac{\tau^{'}}{n}[g(Y,Z)-\{(n-2) M(Y,Z) +mg(Y,Z)\}],
  \end{equation}
 where $S^{"}$ is the Ricci tensor of $\nabla^{"} $ and m is the trace of M(Y,Z).\\
 Take $Y=Z=e_i$ in the above equation we get,
 \begin{equation}
 \tau^{"}=\frac{\tau^{'}}{n}[n-2m(n-1)].
 \end{equation}
 Now from the equation $(7)$, value of $\tau^{'}$ is given as :
 \begin{multline}
 \tau^{'}=c \{\sum_{1\leq i <j\leq n}g(e_j,e_j)g(e_i,e_i)-g(e_i,e_j)g(e_j,e_i)\\ + \sum_{k=1} ^{3} \sum_{1\leq i <j\leq n}(g(\psi_k e_j,e_j)g(\psi_k e_i,e_i)- g(\psi_k e_i,e_j)g(\psi_k e_j,e_i)\\-2g(\psi_k e_i,e_j)g(\psi_k e_j,e_i))\},
\end{multline}
For any p$\in$ $N^{'}$ and X $\in$ $T_p N^{'}$, we have $\psi_i$X $=$ $P_i$ $X$+ $F_i$ $X$ where $P_i$ $\in$ $T_p N^{'}$ and $F_i$ $\in $ $T_p ^{\perp}$.\\
Using this fact in equation $(21)$, we get:
\begin{equation}
 \tau^{'}=c\{n(n-1)+ 3\sum_{k=1} ^3 \|P_k\|^{2}\},
\end{equation}
 where $\|P_k\|^{2}$=$\sum_{i,j=1} ^{n} g^{2}(P_k e_i,e_j)$.\\
 The mean curvature vector H of $N{'}$ in $M^{'}$ is given as:
\begin{equation}
H(p)=\frac{1}{n}\sum_{i=1} ^{n} h(e_i,e_i),
\end{equation}
where h is the second fundamental form of $N^{'}$ in $M^{'}$.\\\\
Also, set;
\begin{equation}
h_{ij} ^{\alpha} =g(h(e_i,e_j),e_{\alpha}),
\end{equation}
where i, j $\in$ $\{1,...,n\}$ and $\alpha$ $\in$ $\{n+1,...,4m\}$.\\\\
The Casorati curvature for the submanifold $N^{'}$ is given as:
\begin{equation}
C=\frac{1}{n}\sum_{\alpha =n+1} ^{4m} \sum_{i,j=1} ^{n} \left(h_{ij} ^{\alpha} \right)^{2}.
\end{equation}
If V is an r-dimensional subspace of $T_p N^{'}$, and let $\{e_1,...,e_r\}$ be an orthonormal basis of V, then we have;
\begin{equation}
C(V)=\frac{1}{r}\sum_{\alpha =n+1} ^{4m} \sum_{i,j=1} ^{r} \left(h_{ij} ^{\alpha} \right)^{2}.
\end{equation} 
The generalized normalized $\delta$-Casorati curvatures $\delta_{C}(r;n-1)$ and $\widehat{\delta_{C}(r;n-1)}$ for any positive real number, r are given as: $[14]$\\\\
$[\delta_{C}(r;n-1)]_p$ $=$ $rC_p$ $+$ $\frac{(n-1)(n+r)(n^{2}-n-r)}{r(n)}$ inf\{C(V)$|$ V is a hyperplane of $T_{p}N^{'},$\}\\\\
if $0<r<n^{2}-n$.\\\\
$[\widehat{\delta_{C}(r;n-1)}]_p$ $=$ $rC_p$ $-$ $\frac{(n-1)(n+r)(r-n^{2}+n)}{r(n)}$ sup\{C(V)$|$ V is a hyperplane of $T_{p}N^{'}$\},\\\\
if $r>n^{2}-n$.\\\\
 The normalized $\delta$-Casorati curvatures $\delta_{C}(n-1)$ and $\widehat{\delta_{C}(n-1)}$ are given as:\\\\
 $[\delta_{C}(n-1)]_p$ $=$ $\frac{1}{2}C_p$ $+$ $\frac{n+1}{2(n)}$ inf\{C(V)$|$ V is a hyperplane of $T_{p}N^{'}$\}.\\\\
 $[\widehat{\delta_{C}(n-1)}]_p$ $=$ $2C_p$ $-$ $\frac{2n-1}{2(n)}$ sup\{C(V)$|$ V is a hyperplane of $T_{p}N^{'}$\}.\\\\
 The above equations imply that normalized $\delta$- Casorati curvature equalities can be obtained from generalized normalized equalities by providing suitable value to the positive number $'r'$ $[26]$.\\\\
 
 Before establishing the main theorems, we will establish the famous lemma of Chen. $[2]$\\\\ 
 \textbf{Lemma 1:} If $n>k\geq 2$ and $a_1,a_2,...a_n,a$ are real numbers such that:
  
  \begin{equation}
  \left(\sum_i ^{n} a_i\right)^2 =(n-k+1)\left(\sum_{i=1} ^n (a_i)^2 +a\right).
  \end{equation}    
  Then,
  \begin{equation}
  2\sum_{1\leq i\leq j \leq k} a_i a_j \geq a,
  \end{equation}
  with equality holding if and only if
   \begin{equation}
  a_1 +a_2 +...+a_k=a_{k+1}=...=a_n.
  \end{equation}
  \section*{Results}
  \section*{(A) Chen invariant, Scalar curvature and Ricci curvature inequalities: }
  \textbf{Theorem 1:} If $N^{'}$ is any submanifold endowed with Ricci quarter-symmetric connection $\nabla$, then $\tau$ follows the inequality as :
  \begin{multline}
  \tau(p) \leq \frac{(n-1)}{2}\left(n\|H\|^{2}+c\{(n-1)+3\sum_{k} ^{3}\frac{\|P_{k}\|^{2}}{n}\}\frac{[n-2m(n-1)]}{(n-1)}\right).
   \end{multline}
   
  \textbf{Proof:}  For $M^{'}$ and $N^{'}$, the Gauss equation is given as:
\begin{equation}
R^{"}(X,Y;Z,W) = R(X,Y;Z,W) +g(h(X,Z),h(Y,W))-g(h(X,W),h(Y,Z)).
\end{equation}
Here '$h$' is the second fundamental form of $N^{'}$. 
Using above equation and equation $(12)$, we get the scalar curvature as:
\begin{equation}
2\tau=c\{(n-1)+3\sum_{k} ^{3} \frac{\|P_{k}\|^{2}}{n}\}[n-2m(n-1)]+n^{2}\|H\|^{2}-\|h\|^{2},
\end{equation}
where $H$=$\frac{1}{n} \sum_{i,j=1} ^n h(e_i,e_j)$ for any orthonormal basis $\{e_1,e_2,...,e_n\}$ of the tangent space $T_p M^{'}$.\\\\
And  $\|h\|^{2} = \sum_{i,j=1} ^{n} g(h(e_i,e_j),h(e_i,e_j)$.\\\\
The equation $(33)$ can be simplified as:
\begin{equation}
n^{2}\|H\|^{2}= 2\tau +c\{(n-1)+3\sum_{k} ^{3}\frac{\|P_{k}\|^{2}}{n}\}[2m(n-1)-n]+\|h\|^{2}. 
\end{equation}
Let $p$ $\in$ $N^{'}$ and ${e_1,...,e_n,...e_m}$ be orthonormal basis at p such that $e_{n+1}$ is parallel to mean curvature vector $H$, then the shape operators are given by:
\begin{equation}
A_{e_{n+1}}=
\begin{pmatrix}
a_1 & 0 & \cdots & o\\
0 & a_2 & \cdots & 0\\
\vdots & \vdots & \ddots & \vdots\\
0 & 0 & \cdots & a_n 
\end{pmatrix}
\end{equation}
$A_{e_r}=h_{ij}^{r}$,   $i,j=1,....,n$; $r=n+2,....,4m$;   $trace(A_r)$ $=$ $0$.\\\\
Considering the equation $(34)$ and the above operators, we get:
 \begin{multline}
 n^{2}\|H\|^{2} = 2\tau(p) +\sum_{i=2} ^n a_i ^{2} +\sum_ {r=n+2} ^{4m} \sum_{i,j=1} ^{n} (h_{ij} ^{r})^{2}\\ + c\{(n-1)+3\sum_{k} ^{3}\frac{\|P_{k}\|^{2}}{n}\}[2m(n-1)-n].
 \end{multline}
 From the algebra, we have:
  \begin{equation}
  \sum_{i<j} (a_i-a_j)^2 =(n-1)\sum_{i} a_i ^2 -2\sum_{i<j} a_i a_j,
  \end{equation}
 \begin{equation}
 \implies n^{2}\|H\|^{2}=\left(\sum_{i} a_i ^2\right)^{2}=\sum_{i} a_i ^2 +2\sum_{i<j}a_ia_j \leq n\sum_{i=1} ^{n} a_i ^{2}.
 \end{equation}
Using above inequality in equation $(36)$, we get:
\begin{multline}
n(n-1)\|H\|^{2} \geq 2\tau(p) + \sum_ {r=n+2} ^{4m} \sum_{i,j=1} ^{n} (h_{ij} ^{r})^{2} +\\+ c\{(n-1)+3\sum_{k} ^{3}\frac{\|P_{k}\|^{2}}{n}\}[2m(n-1)-n],
 \end{multline}
 or, 
\begin{multline}
 \tau(p) \leq \frac{(n-1)}{2}\left(n\|H\|^{2}+c\{(n-1)+3\sum_{k} ^{3}\frac{\|P_{k}\|^{2}}{n}\}\frac{[n-2m(n-1)]}{(n-1)}\right).
\end{multline}
In the above inequality, the equality is established only  when $a_1$=$a_2$=...=$a_n$ and $A_{e_r}$=$0$.\\
$\implies$ p is totally umbilical point.\\\\
\textbf{Remark:} For the real space forms, the part involving $P_k$ will be equal to zero as such the inequality will be same as established in $[11]$.\\\\
 \textbf{Lemma 2}: Let $\zeta $($x_1$, $x_2$,...$x_n$) be a function  in $R^n$ defined by
  \begin{equation}
  \zeta(x_1, x_2,...x_n)=x_1\sum_{i=2} ^n x_i.
  \end{equation}
  If $x_1 +x_2 +...+ x_n= 2\mu$ , then $\zeta$ $\leq$ $\mu^2$.\\
  \textbf{Proof:} From the hypothesis, we have,
  \begin{equation}
  \zeta= x_1(2\mu-x_1),
\end{equation} 
\begin{equation}
\implies\zeta =-(x_1 -\mu)^2 +\mu^2,
\end{equation}
\begin{equation}
\implies\zeta \leq \mu^2.
\end{equation} 
 It follows from the equation $(43)$ that equality is established only when $x_1 = \mu$.\\\\  
 \textbf{Theorem $2$:} If $N^{'}$ is n-dimensional submanifold of m-dimensional Riemannian manifold endowed with Ricci quarter-symmetric connection, then for each unit vector $X$ $\in$ $T_p N^{'}$, the inequality for Ricci curvature is given as: 
 \begin{multline}
Ric(X)\leq c\{(n-1)+3\sum_k ^{3}\sum_{j=2} ^{n} g^{2}(\phi_k X,e_j)\} -c\{(n-1)\\+3\sum_k ^{3} \frac{\|P_{k}\|^{2}}{n} \} [m+(n-2)M(X,X)]+\frac{n^{2}\|H\|^{2}}{4}.
\end{multline}
 \textbf{Proof:} Let $\{e_1,e_2,...,e_n\}$ be the orthonormal basis of $T_p N^{'}$ at any point $p$.\\
 Using our assumption that $N^{'}$ is Einstein manifold and the equations $(32)$ (Gauss Equation)  and $(18)$, we get;
 \begin{multline}
 R(e_1,e_j;e_j,e_1) = R^{*}(e_1,e_j;e_j,e_1)-\frac{\tau{'}}{n} \{M(e_j,e_j)g(e_1,e_1)-M(e_1,e_j)g(e_j,e_1)\\+g(e_j,e_j)M(e_1,e_1)-g(e_1,e_j)M(e_1,e_j)\}+ \sum_{r=n+1} ^{4m}[h_{11} ^rh_{jj} ^{r} - (h_{1j} ^r)^2 ], 
 \end{multline}
 with $g(QX,W)=M(X,W)$.\\\\
If $e_1$=$X$, then the Ricci curvature is defined and given as :
\begin{multline}
Ric(X)=\sum_{i=2}^n R(e_1,e_j,e_j,e_1)=c\{(n-1)+3\sum_k ^{3} g^{2}(\phi_k X,e_j)\} \\-c\{c(n-1)+3\sum_k ^{3} \frac{\|P_{k}\|^{2}}{n} \} [m+(n-2)M(X,X)]\\ +\sum_{r=n+1} ^{4m} \sum_{2} ^{n}[h_{11} ^rh_{jj} ^{r}].
\end{multline}
The above equation implies;
\begin{multline}
Ric(X)\leq c\{(n-1)+3\sum_k ^{3} g^{2}(\phi_k X,e_j)\} -c\{(n-1)\\+3\sum_k ^{3} \frac{\|P_{k}\|^{2}}{n} \} [m+(n-2)M(X,X)]+\sum_{r=n+1} ^{4m} \sum_{2} ^{n}[h_{11} ^rh_{jj} ^{r}]
\end{multline}
Considering the last term of above inequality, we can make a quadratic form given by:
\begin{equation}
\zeta(h_{11}^r,h_{22}^r,...,h_{nn}^r)= \sum_{i=2} ^ n h_{11} ^rh_{ii} ^{r}. 
\end{equation}
Using the Lemma $2$, we have;
\begin{equation}
\zeta (h_{11}^r,h_{22}^r,...,h_{nn}^r) \leq \frac{U_{r} ^2}{4},
\end{equation}
where $U_r = h_{11}^r+h_{22}^r+...+h_{nn}^r$.\\
Using this value and the definition of $H$ in equation $(48)$ we have;
\begin{multline}
Ric(X)\leq c\{(n-1)+3\sum_k ^{3}\sum_{j=2} ^{n} g^{2}(\phi_k X,e_j)\}\\ -c\{(n-1)+3\sum_k ^{3} \frac{\|P_{k}\|^{2}}{n} \} [m+(n-2)M(X,X)]+\frac{n^{2}\|H\|^{2}}{4}.
\end{multline}
 For the condition of equality to hold we should have; 
 \begin{equation}
 h_{1i} ^{r} = 0,
 \end{equation}
 where i $\neq$ $1$ and $\forall$ r.\\
 The another condition follows from the Lemma $2$, which gives,
 \begin{equation}
 h_{11}^r =  h_{22}^r +  h_{33}^r +....+  h_{nn}^r.
 \end{equation}
  If the equality holds for all unit vectors X in $T_p N{'}$, then above inequalities are given as;
  \begin{equation}
  h_{ij} ^{r} = 0, i\neq j,
 \end{equation}
  
  \begin{equation}
  2h_{ii}=h_{11} +h_{22} +....+h_{nn}, i\in \{1,...,n\},
  \end{equation}
  \begin{equation}
  \implies 2h_{11}=2h_{22}=...=2h_{nn}=h_{11} +h_{22} +....+h_{nn}, i\in \{1,...,n\},
  \end{equation}
  
  \begin{equation}
  \implies (n-2)(h_{11}+h_{22}+...+h_{nn})= 0.
   \end{equation}
  For $n \neq 2$,
  \begin{equation}
   h_{ij} ^r =0, i,j=1,2,...,n,r=n+1,...,4m.
  \end{equation}
   Also, 
   \begin{equation}
   h_{ii} ^r=0,
   \end{equation}
   which means the point $p$ is totally geodesic point.\\
  For $n=2$,
  \begin{equation}
  h_{11} ^r =h_{22} ^r, r=3,...,4m,
  \end{equation}
  which means that the point p is totally umbilical.\\\\
\textbf{Remark 1:} If $N(p)$=$\{X\in T_{p} N^{'} :h(X,Y)=0$ $\forall$ $Y \in T_{p} N^{'}\}$, then the above inequality is satisfied by a tangent vector X with $H(p)$ = $0$. \\\\
\textbf{Theorem $3$:} Let $N^{'}$ be n-dimensional submanifold of $4m$-dimensional Quaternionic space form $M^{'}$ with constant sectional curvature given by $4c$. For this submanifold $N^{'}$ which is characterized by Ricci quarter-symmetric metric connection, the inequality between sectional curvature $K$ and scalar curvature $\tau$ at any point p $\in$ $N^{'}$ and section plane $\pi$ is given as:
\begin{multline}
\tau-K(\pi)\leq \frac{n^{2}(n-2)}{2(n-1)} \|H\|^{2}+ \frac{c}{2}\{(n-1)+3\sum_{k} ^{3}\frac{\|P_{k}\|^{2}}{n}\}[n+2mn-2trace(m_{\pi_k
^\perp})] - \\ c\{3\sum_{k=1} ^3 \beta_k (\pi)+1).
\end{multline}

\textbf{Proof:} Let $(e_1,e_2,...,e_n)$ be the orthonormal basis of $T_p N^{'}$ and $(e_{n+1},...,e_m)$ be that of $T_p ^\perp N^{'}$ at any point $p$ such that the section plane $\pi$ is spanned by $(e_1,e_2)$. The mean curvature vector $H$ is taken to be parallel to $e_{n+1}$. Then by the equation of Gauss $(32)$, we have:
\begin{multline}
K(\pi)=K(e_1\wedge e_2)=c(1+3\sum_{k=1} ^3 \beta_k (\pi))+c\{(n-1)+3\sum_{k=1} ^{3} \frac{\|P_k\|^{2}}{n})(trace(m_{\pi_k
^\perp})-m)\\ +h_{11} ^{n+1}h_{22} ^{n+1} +\sum_{r\geq n+2} h_{11} ^{r}h_{22} ^{r}- \left( h_{12} ^{n+1}\right)^{2} -\sum_{r\geq n+2} \left( h_{12}^{r}\right)^{2},
\end{multline} 
where m is the trace of $M$, $trace(m_{\pi_k ^{\perp}})$ is equal to ($m- M(e_1,e_1)-M(e_2,e_2))$ and $\beta_k(\pi)$=$\sum_{1\leq i \leq j \leq 2} g^{2}(P_k e_i,e_j)$ with $k=1,2,3$.\\\\
If we assume;
\begin{multline}
\mu= 2\tau-c((n-1)+3\sum_{k} ^{3}\frac{\|P_{k}\|^{2}}{n})[n-2m(n-1)]-\frac{n^{2}(n-2)}{n-1} \|H\|^{2}.
\end{multline}
Then equation $(33)$ and above equation gives;
\begin{equation}
n^{2}\|H\|^{2}=(n-1)(\mu+\|h\|^{2}).
\end{equation}
The above implies;
\begin{multline}
\left(\sum_{i=1} ^n h_{ii} ^{n+1}\right)^2 =(n-1)\left(\sum_{i=1} ^n (h_{ii} ^{n+1})^2 +\sum_{r=n+2}^{4m} \sum_{i,j=1} ^n (h_{ij} ^ r)^2 +\mu \right).
 \end{multline}
 Using Lemma $1$,we get;
   \begin{equation}
   2\sum_{1\leq i \leq j \leq k} h_{ii} ^{n+1} hh_{jj} ^{n+1} \geq \mu + \sum_{i \neq j} (h_{ij} ^{n+1})^{2}+ \sum_{r=n+2}^m \sum_{i,j=1} ^n (h_{ij} ^ r)^2.
   \end{equation}
   Using these values in the equation involving  $K(\pi)$, we get;
 \begin{multline}
   K(\pi)\geq c(1+3\sum_{k=1} ^3 \beta_k (\pi))+c\{(n-1)+3\sum_{k=1} ^{3} \frac{\|P_k\|^{2}}{n}\}(trace(m_{\pi_k
^\perp})-m)+\frac{\mu}{2} \\+ \sum_{r=n+1}^{4m} \sum_{j>2} ( (h_{1j} ^ r)^2+(h_{2j} ^ r)^2)+\frac{1}{2}\sum_{i\neq j>2} ( (h_{ij} ^ {n+1})^2+\frac{1}{2}\sum_{r=n+2} ^{4m}\sum_{i,j>2} ( (h_{ij} ^ {r})^2\\+\frac{1}{2}\sum_{r=n+2} ^{4m} (h_{11} ^ {r} + h_{22} ^ {r})^{2}.
 \end{multline}
 On solving above inequality we get,
 \begin{multline}
 \tau-K(\pi)\leq \frac{n^{2}(n-2)}{2(n-1)} \|H\|^{2}+ \frac{c}{2}\{(n-1)+3\sum_{k} ^{3}\frac{\|P_{k}\|^{2}}{n}\}[n-2m(n-1)]\\ -  c\{ 1+3\sum_{k=1} ^3 \beta_k (\pi))-c\{3\sum_{k=1} ^{3} \frac{\|P_k\|^{2}}{n}+(n-1)\}(trace(m_{\pi_k
^\perp})-m).
 \end{multline} 
 In the above inequality, the equality condition holds if the terms involving $h_{ij} ^{r}$ is equal to zero and as such we have;\\\\
 $$h_{1j} ^{n+1}=h_{2j} ^{n+1}=0,$$
 $$h_{ij} ^{n+1} =0; i\neq
 j>2,$$
 $$h_{1j} ^{r}=h_{2j} ^{r}=...=h_{ij} ^{r}; r=n+2,....4m;i,j>3,$$
 $$h_{11} ^{r} +h_{22} ^{r} =0;r=n+2,...,4m.$$
 Choose $e_1$ and $e_2$ such that $h_{12} ^{n+1} =0$. Also by Chen's Lemma $1$, we have;
 $$h_{11} ^{n+1} +h_{22} ^{n+1}=h_{33} ^{n+1}=...=h_{nn} ^{n+1}.$$
 The above equations imply that the shape operator take the form as :
 \begin{equation}
    A_{e_{n+1}}=
\begin{pmatrix}
T & 0 & \cdots & 0 & 0\\
0 & U & \cdots & 0 & 0\\
0 & 0 & I & \cdots & 0\\
\vdots & \vdots & \ddots & \vdots\\
0 & 0 & \cdots & 0 &I 
\end{pmatrix}
    \end{equation}
     \begin{equation}
    A_{e_{r}}=
\begin{pmatrix}
T_r &U_r  & 0\\
U_r & -U_r & 0\\
0 & 0 & 0_{4m-2}\\
\end{pmatrix}
    \end{equation}
    where $T+U=I$ and $T_r$, $U_r$ $\in$ R.\\\\
   
 \textbf{Theorem $4$:} If $N^{'}$ is n-dimensional submanifold of $4$m-dimensional quaternionic projective space $M^{'}$ ($c>0$) with constant quaternionic sectional curvature given by $4$c and endowed with Ricci quarter-symmetric connection, then for any plane $\pi$ $\in$ $T_p N^{'}$, the inequality for Chen's invariant $\delta$ is given as:\\
 \begin{equation}
 \delta_{N^{'}}\leq \frac{n^{2}(n-2)}{2(n-1)} \|H\|^{2}
+ \frac{c}{2}(n+8)(n+2mn-2trace(m_{\pi_k
^\perp}))-c,
 \end{equation}
 with the equality holding if and only if $N^{'}$ is invariant.\\\\
 \textbf{Proof:} For the positive sectional curvature($c>0$), we have to maximize the term $\frac{3c(n-1)}{2}\sum_{k=1}^{3} \frac{\|P_k\|^{2}}{n}(n+2mn-2trace(m_{\pi_k
^\perp}))-\frac{3c}{2}\sum_{k=1} ^{3} \beta
_k(\pi)$ of inequality $(61)$.\\
For constant c,m and $trace(m_{\pi_k
^\perp})$, the maximum value can be achieved if $\|P_k\|^{2}$ =$n$ and $ \beta
_k(\pi)=0$.\\\\
Using these values, we get the required inequality.\\\\
 \textbf{Remark 2:}The condition,$\|P_k\|^{2}$ =$n$ and  $ \beta
_k(\pi)=0$ implies $N^{'}$ is invariant.\\\\
\textbf{Remark 3:} If $M(X,Y)=0$, we get the inequality as in $[12]$.\\\\
\textbf{Theorem $5$:} If $N^{'}$ is n-dimensional submanifold of $4$m-dimensional quaternionic hyperbolic space $M^{'}$ ($c<0$) with constant quaternionic sectional curvature given by $4$c and endowed with Ricci quarter-symmetric connection, then for any plane $\pi$ $\in$ $T_p N^{'}$, the inequality for Chen's invariant $\delta$ is given as:\\
 \begin{equation}
 \delta_{N^{'}}\leq \frac{n^{2}(n-2)}{2(n-1)} \|H\|^{2}
+ \frac{c}{2}(n-1)(n+2mn-2trace(m_{\pi_k
^\perp}))-c,
 \end{equation}
 With the equality holding if and only if $N^{'}$ is anti quasi-invariant.\\\\
 \textbf{Proof:} For $c<0$, the estimate for $\delta_{N^{'}}$ can be obtained by minimizing the term $\frac{3c(n-1)}{2}\sum_{k=1}^{3} \frac{\|P_k\|^{2}}{n}(n+2mn-2trace(m_{\pi_k
^\perp}))-\frac{3c}{2}\sum_{k=1} ^{3} \beta
_k(\pi)$ of inequality $(61)$.\\\\
By simplification it follows that $\|P_k\|^{2}$ =$0$ and $ \beta
_k(\pi)=0$.\\\\
The above conclusion exhibits that span of $\pi$ $=$ $span(e_1,e_2)$ is orthogonal to span of $\{span(\phi_{k}e_{i}$ $|=3,...,n$\}
\section*{(B) Casorati curvature Inequalities:}
\textbf{Theorem 1:} If $N^{'}$ is $n$-dimensional submanifold in $4m$-dimensional quaternionic Kaehler manifold $M^{'}$. Then:\\\\
(i). The generalized normalized $\delta$-Casorati curvature $\delta_C(r,n-1)$ satisfies;
  \begin{multline}
 \delta_C(r,n-1) \geq n(n-1)\rho\\ - c\{(n-1)+3\sum_{k} ^{3} \frac{\|P_{k}\|^{2}}{n}\}[n-2m(n-1)].
  \end{multline}
  (ii).  The generalized normalized $\delta$-Casorati curvature $\widehat{\delta_C (r,n-1)}$ satisfies;
  \begin{multline}
 \widehat{\delta_C (r,m)}  \geq n(n-1)\rho\\ - c\{(n-1)+3\sum_{k} ^{3} \frac{\|P_{k}\|^{2}}{n}\}[n-2m(n-1)].
  \end{multline}
  In the inequality $(73)$ , the equality holds if and only if the submanifold M is invariantly quasi-umbilical with trivial normal connection in $M^{'}$, such that the shape operator $A_{e_i}$ with $i$ $\in$ $\{n+1,...,4m\}$ with respect to suitable orthonormal tangent frame $\{e_1,...,e_{n}\}$ and normal orthonormal frame $\{e_{n+1},...,e_{4m}\}$ takes the form as:
  \begin{equation}
  A_{e_{n+1}}=
\begin{pmatrix}
u & 0 & 0 & \cdots & 0 & 0\\
0 & u & 0 & \cdots & 0 & 0\\
0 & 0 & u & \cdots & 0 & 0\\
\vdots & \vdots & \ddots & \vdots\\
0 & 0 & 0 & \cdots & u & 0\\
0 & 0 & 0 & \cdots & 0 &  \frac{n(n-1)}{r}u
\end{pmatrix}
  \end{equation}
  \begin{equation}
  A_{e_{n+2}}=...,=A_{e_{4m}}=0.
  \end{equation}
   
\textbf{Proof:} By equation $(33)$ we have;
\begin{equation}
2\tau=c\{(n-1)+3\sum_{k} ^{3} \frac{\|P_{k}\|^{2}}{n}\}[n-2m(n-1)]+n^{2}\|H\|^{2}-\|h\|^{2},
\end{equation}
where $H$=$\frac{1}{n} \sum_{i,j=1} ^n h(e_i,e_j)$ for any orthonormal basis $\{e_1,e_2,...,e_n\}$ of the tangent space $T_p M^{'}$ and m is the trace of $M(X,Y)$.\\\\
Now, consider the function T which is associated with the following polynomial in the components $(h_{ij} ^{\alpha}), i,j=1,...,n;\alpha=n+1,...,4m$ of second fundamental form h of $N^{'}$ in $M^{'}$:\\
\begin{multline}
T= rC +\frac{(n-1)(n+r)(n^{2}-n-r)}{rn}C(L)-2\tau\\ + c\{(n-1)+3\sum_{k} ^{3} \frac{\|P_{k}\|^{2}}{n}\}[n-2m(n-1)],
\end{multline}
where L is the hyperplane of $T_{p}N^{'}$.\\\\
If L is assumed to be spanned by $e_1,...,e_{n-1}$, then we have:
\begin{multline}
T=\frac{r}{n}\sum_{\alpha=n+1} ^{4m} \sum_{i,j=1} ^{n} (h_{ij})^{2} +\frac{(n+r)(n^{2}-n-r)}{rn}\sum_{\alpha=n+1} ^{4m} \sum_{i,j=1} ^{n} (h_{ij})^{2} -2\tau\\ + c\{(n-1)+3\sum_{k} ^{3} \frac{\|P_{k}\|^{2}}{n}\}[n-2m(n-1)].
\end{multline}
Using equation $(77)$ we get;
\begin{multline}
T=\frac{n+r}{n}\sum_{\alpha=n+1} ^{4m} \sum_{i,j=1} ^{n} (h_{ij})^{2}\\ +\frac{(n+r)(n^{2}-n-r)}{rn}\sum_{\alpha=n+1} ^{4m} \sum_{i,j=1} ^{n} (h_{ij})^{2}-\sum_{\alpha=n+1} ^{4m}\left(\sum_{i+1} ^{n}h_{ii}^{\alpha}\right)^{2}.
\end{multline}
The above equation can be written as:
\begin{multline}
T=\sum_{\alpha=n+1} ^{4m}\sum_{i=1} ^{n-1}\left[\frac{(n-1)^{2}+(n-1)(r+1)-r)}{r} (h_{ii}^{\alpha})^{2} +\frac{2(n+r)}{n}(h_{i(n-1)+1}^{\alpha})^{2}\right]\\ \sum_{\alpha=n+1}^{4m}\left[\frac{2(n-1)(n+r)}{r}\sum_{1 \leq i<j\leq m-1} (h_{ij} ^{\alpha})^{2}-2 \sum_{1 \leq i<j\leq n} h_{ii} ^{\alpha} h_{jj} ^{\alpha} \right]\\
+\frac{r}{n}\sum_{\alpha=n+1} ^{4m} (h_{n n} ^{\alpha})^{2}.
\end{multline}
Hence T is a quadratic polynomial in the components of the second fundamental form and from $(81)$, we deduce the critical points;\\\\
$h^{c}$= $\left[ h_{11}^{n},h_{12}^{n+1},...,h_{n n}^{n+1},...,h_{11}^{4m},h_{12}^{4m},...,h_{n n}^{4m} \right].$\\\\
of T are the solutions of the following system of linear homogeneous equations:
\begin{equation}
\frac{\partial T}{\partial h_{ii}^{\alpha}}=\frac{2(n-1)(n+r)}{r}h_{ii}^{\alpha}-2\sum_{k=1}^{n}h_{kk}^{\alpha}=0,
\end{equation}
\begin{equation}
\frac{\partial T}{\partial h_{n n}^{\alpha}}=\frac{2r}{n}h_{n}^{\alpha}-2\sum_{k=1}^{n}h_{kk}^{\alpha}=0,
\end{equation}
\begin{equation}
\frac{\partial T}{\partial h_{ij}^{\alpha}}=\frac{4(n-1)(n+r)}{r}h_{ij}^{\alpha}=0,
\end{equation}
\begin{equation}
\frac{\partial T}{\partial h_{i n}^{\alpha}}=\frac{4(n+r)}{n}h_{i n}^{\alpha}=0,
\end{equation}
with $i,j \in \{1,...,n-1\}$, $i\neq j$ and $\alpha \in \{n+1,...,4m\}$.\\\\
The above system of equations shows that every solution $h^{c}$ has $h_{ij}^{\alpha}=0$ for $i \neq j$. The Hessian matrix can be computed as:
\begin{equation}
  H(p)=
\begin{pmatrix}
H_1 & 0 & 0\\
0 & H_2 & 0\\
0 & 0 & H_3
\end{pmatrix}
  \end{equation}
  Here $H_1$ is given as:
  \begin{equation}
  H_1=
\begin{pmatrix}
\frac{2(n-1)(n+r)}{r} -2 & -2 & -2 & \cdots & -2 & -2\\
-2 & \frac{2(n-1)(n+r)}{r} -2 & -2 & \cdots & -2 & -2\\
-2 & -2 & \frac{2(n-1)(n+r)}{r} -2 & \cdots & -2 & -2\\
\vdots & \vdots & \ddots & \vdots\\
-2 & -2 & -2 & \cdots & \frac{2(n-1)(n+r)}{r} -2 & -2\\
-2 & -2 & -2 & \cdots & -2 &  \frac{2r}{n}
\end{pmatrix}
  \end{equation}
  $H_2$ $=$ $diag\left(\frac{4(n-1)(n+r)}{r},\frac{4(n-1)(n+r)}{r},...,\frac{4(n-1)(n+r)}{r}\right)$,\\\\
  $H_3$ $=$ $diag\left(\frac{4(n+r)}{n},\frac{4(n+r)}{n},...,\frac{4(n+r)}{n}\right)$.\\\\
 It can be found by direct computation that the Hessian matrix $H(T)$ of T has the eigenvalues as follows:\\\\
$\lambda_{11}=0$,  $\lambda_{22}=\frac{2(n-1)n^{2} +r^{2})}{rn}$ , $\lambda_{33}=...,=\lambda_{m+1 m+1}=\frac{2(n-1)(n+r)}{r}$,\\\\
$\lambda_{ij} =\frac{4(n-1)(n+r)}{r}$,$\lambda_{i(n-1)+1} =\frac{4(n+r)}{n}$;, $\forall$ i,j $\in$ $\{1,...,n-1\}$, $i\neq j $.\\\\
It follows that the Hessian matrix is positive semi-definite and admits one eigenvalue equal to zero. Hence, we can see that T is parabolic and reaches a minimum at $h^{c}$. In fact, at the critical point $h^{c}$, the function T reaches the global minimum. So, we have $T(h^{c})=0$. Therefore, we deduce $T\geq 0$ and this implies;
\begin{multline}
2\tau \leq rC +\frac{(n-1)(n+r)(n^{2}-n-r)}{rn}C(L)\\+ c\{(n-1)+3\sum_{k} ^{3} \frac{\|P_{k}\|^{2}}{n}\}[n-2m(n-1)].
\end{multline}
We know that normalized scalar curvature $\rho$ is given as:
\begin{equation}
\rho=\frac{2\tau}{n(n-1)}.
\end{equation}
Using this value in above equation we get;
\begin{multline}
\rho \leq \frac{r}{n(n-1)} C +\frac{(n+r)(n^{2}-m-r)}{rn^{2}}C(L)\\ + c\{(n-1)+3\sum_{k} ^{3} \frac{\|P_{k}\|^{2}}{n}\} \frac{[n-2m(n-1)]}{n(n-1)}.
\end{multline}
The above statement is equivalent to;
\begin{multline}
rC +\frac{(n-1)(n+r)(n^{2}-m-r)}{rn}C(L) \geq n(n-1)\rho\\ - c\{(n-1)+3\sum_{k} ^{3} \frac{\|P_{k}\|^{2}}{n}\}[n-2m(n-1)].
\end{multline}
The required inequalities follow by taking the infimum and supremum respectively, over all the tangent hyperplanes $T_{p}N^{'}$.\\\\ On the other hand, it is clear that the equality sign follows if and only if;
$h_{ij}^{\alpha}=0$, $\forall$ i,j $\in$ $\{1,...,n\}$, $i \neq j$:\\
and
\begin{equation}
h_{n n}^{\alpha} = \frac{n(n-1)}{r}h_{11}^{\alpha}=,...,\frac{n(n-1)} {r}h_{n-1 n-1}^{\alpha},
\end{equation} 
for all $\alpha$ $\in$ $\{n+1,...,4m\}$.\\\\
  Since, $h(e_i, \zeta_p)$ $\neq$ $0$, we have that the equality holds only when the submanifold  is invariantly quasi-umbilical and shape operator takes the form given in equations $(75)$ and $(76)$.\\ 
  Similarly, inequality $(74)$ can be proved.\\\\ 
 \textbf{Corollary $1.1$:} If $N^{'}$ is $n$-dimensional submanifold in $4m$-dimensional quaternionic Kaehler manifold. Then:\\\\
(i). The normalized $\delta$-Casorati curvature $\delta_C(n-1)$ satisfies;
  \begin{multline}
 \delta_C(n-1) \geq \rho - c\{(n-1)+3\sum_{k} ^{3} \frac{\|P_{k}\|^{2}}{n}\}\frac{[n-2m(n-1)]}{n(n-1)}.
  \end{multline}
(ii). The normalized $\delta$-Casorati curvature $\widehat{\delta_C (n-1)}$ satisfies;
  \begin{multline}
 \widehat{\delta_C (n-1)} \geq \rho - c\{(n-1)+3\sum_{k} ^{3} \frac{\|P_{k}\|^{2}}{n}\}\frac{[n-2m(n-1)]}{n(n-1)}.
  \end{multline}
  In the inequalities $(93)$ and $(94)$ , the equality holds if and only if the submanifold M is invariantly quasi-umbilical with trivial normal connection in $M^{'}$, such that the shape operator $A_{e_i}$ with $i$ $\in$ $\{n+1,...,4m\}$ with respect to suitable orthonormal tangent frame $\{e_1,...,e_{n}\}$ and normal orthonormal frame $\{e_{n+1},...,e_{4m}\}$ takes the form as:
  \begin{equation}
  A_{e_{n+1}}=
\begin{pmatrix}
u & 0 & 0 & \cdots & 0 & 0\\
0 & u & 0 & \cdots & 0 & 0\\
0 & 0 & u & \cdots & 0 & 0\\
\vdots & \vdots & \ddots & \vdots\\
0 & 0 & 0 & \cdots & u & 0\\
0 & 0 & 0 & \cdots & 0 &  2u
\end{pmatrix}
  \end{equation}
  \begin{equation}
  A_{e_{n+2}}=...,=A_{e_{4m}}=0.
  \end{equation}
  and \begin{equation}
  A_{e_{n+1}}=
\begin{pmatrix}
2u & 0 & 0 & \cdots & 0 & 0\\
0 & 2u & 0 & \cdots & 0 & 0\\
0 & 0 & 2u & \cdots & 0 & 0\\
\vdots & \vdots & \ddots & \vdots\\
0 & 0 & 0 & \cdots & 2u & 0\\
0 & 0 & 0 & \cdots & 0 &  u
\end{pmatrix}
  \end{equation}
  \begin{equation}
  A_{e_{n+2}}=...,=A_{e_{4m}}=0
  \end{equation}
  \textbf{Proof:} 
  It can be easily seen that:
  \begin{equation}
  [\delta_C(\frac{n(n-1)}{2},m)]_p=n(n-1)[\delta_C(n-1)]_p,
  \end{equation}
  And, 
  \begin{equation}
  [\widehat{\delta_C (2n(n-1),m)}]_p=n(n-1)[\widehat{\delta_C (n-1)}]_p.
  \end{equation}
  The proof follows by putting $r=\frac{n(n-1)}{2}$ in $(73)$ and $r=2n(n-1)$ in $(74)$ and using above two equations.\\\\

\end{document}